\documentclass[12pt]{article}
\usepackage{float}

\usepackage{amsmath, amssymb}
\usepackage{geometry}
\usepackage{enumitem}
\usepackage{hyperref}
\usepackage{titlesec}
\usepackage{graphicx}
\usepackage{cite}
\usepackage{booktabs}
\usepackage{epstopdf} 
\usepackage{version}
\usepackage{bm}
\usepackage{enumerate}
\usepackage{epsfig}
\usepackage{caption}
\usepackage{xcolor}
\usepackage{setspace}
\usepackage{caption}
\usepackage{float}
\DeclareCaptionLabelFormat{cont}{#1~#2\alph{ContinuedFloat}}
\usepackage{footnote}
\newtheorem{theorem}{Theorem}[section]

\newtheorem{lemma}[theorem]{Lemma}

\newtheorem{definition}[theorem]{Definition}
\newtheorem{example}[theorem]{Example}

\newtheorem{remark}[theorem]{Remark}

\numberwithin{equation}{section}

\excludeversion{comment1}

\titleformat{\section}{\large\bfseries}{\thesection.}{0.5em}{}

\usepackage{float}
\begin{document}
\title{\textbf{Accelerating Convergence in Series and Infinite Integrals: \\Revisiting Levin and Sidi's Contributions}}

\author{ David Levin }
\author{David Levin\\
        School of Mathematical Sciences, Tel-Aviv University, Israel}
\date{}

\date{}
\maketitle

\begin{abstract}
The evaluation of slowly converging series and infinite integrals is a key challenge in numerical analysis and computational mathematics. In their influential 1981 paper, the author and Avram Sidi introduced two effective nonlinear transformations, the d-transformation for series and the D-transformation for infinite integrals, aimed at speeding up their convergence. This review summarizes, contextualizes, and evaluates their contributions, highlighting the mathematical basis, practical significance, and legacy of their work.
\end{abstract}

\section{Introduction}

Infinite series and integrals arise frequently in applied mathematics, physics, and engineering, but their numerical evaluation is often hampered by slow convergence or even divergence due to insufficient decay of the terms. While classical acceleration techniques can be effective in specific instances, they often lack the generality needed for broader applicability. In a seminal 1981 work, Levin and Sidi introduced two nonlinear transformations that exploit the asymptotic behavior of summands and integrands. These transformations yield significant convergence acceleration across a wide range of problems commonly encountered in applied analysis and computational science.

This paper revisits their contribution, situating it within the broader context of summation and integration methods. It highlights both the mathematical innovation of their approach and its practical impact on numerical computation.

\section{Background and Classical Methods}

Brezinski’s book \cite{Brezinski} serves as a comprehensive reference on convergence acceleration techniques, offering extensive historical insight and broad methodological coverage. 
In contrast, Sidi’s book \cite{Sidi2003} presents a focused and systematic treatment of asymptotics-based methods for accelerating the convergence of series and integrals, with particular emphasis on their theoretical foundations and practical implementation.

In this paper, we investigate the two closely related problems: the evaluation of infinite series
\begin{equation}\label{S}
S = \sum_{n=0}^\infty f(n),
\end{equation}
and the evaluation of infinite integrals
\begin{equation}\label{I}
I = \int_{0}^\infty f(t)dt.
\end{equation}

Several classical techniques have been developed to accelerate the convergence of series:

\begin{itemize}
    \item Euler transformation (for alternating series \cite{Euler})
    \item Aitken's $\Delta^2$ process \cite{Aitken}
    \item Padé approximants \cite{Pade1892}
    \item Wynn’s $\varepsilon$-algorithm \cite{Wynn}
    \item Levin-type $t$- and $u$-transformations \cite{Levin1972}
    \item Gray, Atchison and McWilliams' G-transformation (For infinite integrals \cite{Gray}).
\end{itemize}

These methods are typically linear or semi-linear and often rely on specific properties of the function 
$f$, such as alternating behavior or slow decay. However, their effectiveness can diminish when applied to more complex functions, particularly when asymptotic information is not explicitly taken into account.

A Levin-type transformation is a non-linear method designed to improve the convergence of series or approximate their sums (even when divergent) by using explicit knowledge or assumptions about the asymptotic behavior of the remainder terms. Specifically, we are concerned with the behavior, as $N\to \infty$, of the remainder terms
\begin{equation}\label{RS}
R_S(N) = \sum_{n=N}^\infty f(n)=S-\sum_{n=0}^{N-1} f(n),
\end{equation}
in the case of a series, and
\begin{equation}\label{RI}
R_I(x) = \int_{x}^\infty f(t)dt=I-\int_{0}^x f(t)dt,
\end{equation}
in the case of infinite integrals. Levin-type transformations are particularly well-suited for cases where the remainder terms admit simple asymptotic expansions. In the following, we explore the connection between asymptotic models that describe the behavior of 
$f$ for large arguments and their implications for analyzing the corresponding remainder terms. We begin with the case of infinite series and subsequently extend the discussion to infinite integrals.

\begin{subsection}{A linear model with constant coefficients}

 Consider the infinite series \eqref{S}. We begin by assuming that the function \( f \) satisfies a linear recurrence relation of order \( k \) with constant coefficients:
\begin{equation} \label{linearconstant}
a_0 f(n) + a_1 f(n+1) + \cdots + a_k f(n+k) = 0, \quad n \ge 0,
\end{equation}
where \( a_0, a_1, \dots, a_k \) are constants. It is worth noting that both Padé approximants~\cite{Pade1892} and Wynn’s \( \varepsilon \)-algorithm~\cite{Wynn} are capable of recovering the exact sum of an infinite series when the summand satisfies a relation of the form \eqref{linearconstant}. In what follows, we present a simple and direct method for computing the exact sum using the remainder framework.

The general solution to \eqref{linearconstant} is given by
\begin{equation}
f(n) = \sum_{j=1}^m p_j(n) r_j^n,
\end{equation}
where each \( p_j(n) \) is a polynomial of degree \( m_j - 1 \), where $\sum_{j=1}^m m_j = k$. If \( r_j \) is complex, the terms \( r_j^n \) and \( \overline{r_j}^n \) can be combined into real-valued trigonometric functions.

Thus, the class of functions includes:
\begin{itemize}
  \item Geometric sequences: \( r^n \),
  \item Polynomial sequences,
  \item Polynomially weighted exponentials: \( n^k r^n \),
  \item Trigonometric sequences: \( \cos(\omega n), \sin(\omega n) \).
\end{itemize}

Considering \eqref{RS}, we apply the summation from \( n = N \) to \( \infty \) to the recurrence relation \eqref{linearconstant}, yielding a linear recurrence relation for the remainders:
\begin{equation}\label{summing}
a_0R_s(N)+a_1R_S(N+1)+\cdots+a_kR_S(N+k)=0.
\end{equation}
Let \( A_M = \sum_{n=0}^{M-1} f(n) \) denote the partial sum up to index \( M-1 \). Recalling the definition of the remainder term \eqref{RS}, we can write
\begin{equation} \label{summing2}
\sum_{j=0}^k a_j \left(S - A_{N+j} \right) = 0.
\end{equation}
Solving for \( S \), we obtain the expression for the infinite sum \eqref{S}, assuming that the sequence \( \{f(n)\} \) satisfies the recurrence relation \eqref{linearconstant}:
\begin{equation} \label{sumf}
S = \frac{\sum_{j=0}^k a_j A_{N+j}}{\sum_{j=0}^k a_j}.
\end{equation}

The above example demonstrates how the information about the behavior of the remainder $R_S(N)$ can be effectively used to estimate the value of the infinite sum $S$. 
    
\end{subsection}
\section{The Work of Levin and Sidi}

Levin and Sidi introduced:
\begin{itemize}
\item The \textbf{d-transformation} for infinite series
\item The \textbf{D-transformation} for infinite integrals
\end{itemize}

These transformations are based on the assumption that the function \( f \) satisfies either a linear difference equation (in the case of series) or a linear differential equation (in the case of integrals), with coefficients that admit specific asymptotic expansions.

A key advantage of this approach is that explicit knowledge  
of the governing equation is not required. It is enough to assume  
its existence and estimate its order. This makes the transformations  
widely applicable.

The construction of these transformations relies on deriving  
general asymptotic expansions for the remainder terms in series  
and integrals, based on the assumed behavior of the function.  
In the following, we outline the main steps in the derivation  
of the Levin and Sidi transformations.

\subsection{Mathematical Foundations}

As noted above, the \( d \)-transformation is developed under the assumption  
that the sequence \( \{f(n)\} \) satisfies a linear difference equation whose coefficients  
belong to a specific class of asymptotic functions, defined as follows:

\begin{definition}{\bf Asymptotic class $A^{(\gamma)}$}

Let \( \gamma \in \mathbb{R} \). A function \( p \) is said to belong to the class \( A^{(\gamma)} \) if it admits an asymptotic expansion of the form
\begin{equation}
p(x) \sim \sum_{i=0}^{\infty} \alpha_i x^{\gamma - i}, \quad \text{as } x \to \infty,
\end{equation}
where the \( \alpha_i \) are constants.

This means that for each \( n \geq 0 \), there exists a constant \( C_n \) such that
\[
\big| p(x) - \sum_{i=0}^{n} \alpha_i x^{\gamma - i}\big| \leq C_n x^{\gamma - n - 1}, \quad \text{for all sufficiently large } x.
\]
\end{definition}

Specifically, the \( d \)-transformation is particularly effective for accelerating the convergence of sequences whose elements belong to the classes of functions \( \{B^{(m)}\} \), defined as follows:

\begin{definition}{\bf The set $B^{(m)}$}

Let \( m \in \mathbb{N} \). The set \( B^{(m)} \) consists of all infinite sequences \( \{f(n)\}_{n=0}^{\infty} \) whose elements satisfy a linear difference equation of order \( m \), that is,
\begin{equation}
f(n) = \sum_{k=1}^{m} p_k(n) \Delta^k f(n),
\end{equation}
where \( \Delta f(n) = f(n+1) - f(n) \), and the coefficient functions \( p_k(n) \) are such that, when considered as functions of a continuous variable, belong to the asymptotic class \( A^{(k)} \), for \( k = 1, 2, \dots, m \).
\end{definition}
\begin{remark}
While linear models with constant coefficients are  
limited to sums of exponentials, polynomials, and  
trigonometric functions, the sets \(\{B^{(m)}\}\)  
include a much broader class of sequences, such as  
Bessel functions, orthogonal polynomials, and  
hypergeometric functions.
\end{remark}

The following lemma, which is proved in \cite{LevinSidi1981}, assists in identifying the set $B^{(m)}$ to which certain sequences of terms of given infinite series belong.

\vspace{1em}

\begin{lemma}\label{LemmaB}{\bf Algebraic rules for the sets \(\{B^{(m)}\)\} }

If $\{f(n)\}\in B^{(m)}$ and $\{g(n)\}\in B^{(k)}$ then
$\{f(n+g(n)\}\in B^{(m+k)}$,
$\{f(n)g(n)\}\in B^{(mk)}$ and $\{f^2(n)\}\in B^{(\frac{m(m+1)}{2})}$.
\end{lemma}

The following theorem, proved in \cite{LevinSidi1981}, is the principal result underlying the definition of the 
$d$-transformation. It establishes an asymptotic expansion for the remainder of an infinite series whose term sequences belong to the class  $B^{(m)}$.

\begin{theorem}\label{TRS}[Theorem 2 in \cite{LevinSidi1981}]

Let $|\sum_{n=0}^\infty f(n)|<\infty$, and let $\{f(n)\}\in B^{(m)}$. If 
\begin{equation}
\lim_{n\to\infty}[\Delta^ip_k(n)][\Delta^{k-i}f(n)]=0, \ \ k=i,i+1,\cdots,m,\ \ i=1,2,\cdots,m,
\end{equation}
and 
\begin{equation}
\sum_{k=1}^m\ell(\ell-1)\cdots(\ell-k+1)p_{k,0}\ne 1,\ \ell\ge -1, \ \ell\ integer,
\end{equation}
where $p_{k,0}=\lim_{n\to\infty}n^{-k}p_k(n)$, then the remainder $R_S(N)=\sum_N^\infty f(n)$, has an asymptotic expansion of the form
\begin{equation}\label{Rem}
R_S(N)\sim \sum_{k=0}^{m-1}\Delta^kf(N)N^k\left(\beta_{k,0}+\frac{\beta_{k,1}}{N}+\frac{\beta_{k,1}}{N^2}+\cdots \right), \ \ {\text as}\ N\to \infty.
\end{equation}
\end{theorem}

\subsection{The \texorpdfstring{$d$}{d}-transformation for infinite series}

We now utilize the asymptotic expansion~\eqref{Rem} to derive the \( d \)-transformation for infinite series $S = \sum_{n=0}^\infty f(n)$. This is achieved by truncating each of the \( m \) infinite sums appearing in~\eqref{Rem} after \( r \) terms. Recalling that \( R_S(N) = S - \sum_{n=0}^{N-1} f(n) \) and denoting the partial sums of the series by $A_N(f)=\sum_{n=0}^{N-1} f(n)$, we obtain, up to the truncation error,
\begin{equation} \label{RemT}
S \approx A_N(f)+ \sum_{k=0}^{m-1} \Delta^k f(N) \, N^k \sum_{i=0}^{r-1} \frac{\beta_{k,i}}{N^i}.
\end{equation}
There are \( mr + 1 \) unknown parameters in~\eqref{RemT}, namely the sum \( S \) and the coefficients \( \{ \beta_{k,i} \} \).
To determine these parameters, we construct a system of \( mr + 1 \) equations of the form
\begin{equation}
\label{RemTNj}
\bar S = A_{N_j}(f) + \sum_{k=0}^{m-1} \Delta^k f(N_j) \, N_j^k \sum_{i=0}^{r-1} \frac{\bar{\beta}_{k,i}}{N_j^i}, \quad j = 1,2,\ldots, mr+1,
\end{equation}
where \( 0 < N_1 < N_2 < \cdots < N_{mr+1} \). A common choice is \( N_j = \ell + j \),  
for \( j = 1,2,\ldots, mr+1 \).  
With this choice, we solve system~\eqref{RemTNj},
and denote the result for \( \bar{S} \) by \( d^{(m)}_{r,\ell} \).  

This quantity is known as the \( d \)-\emph{transformation}  
of the sequence of partial sums of the series,  
and is denoted by \( d^{(m)}_{r,\ell}(\{A_N(f)\}) \).  
The parameters involved are:  
\begin{itemize}
  \item \( m \): the order of the difference operator;  
  \item \( r \): the number of terms used in the asymptotic expansions;  
  \item \( \ell \): an index determining the sequence \( \{N_j\} \).  
\end{itemize}

Note that computing \( d^{(m)}_{r,\ell} \)  
requires \( \ell + mr + 1 \) terms from the series.  
For the special case \( \ell = 0 \),  
we write the transformation as  
\( d^{(m)}_{r}(\{A_N(f)\}) \).

\subsection{Examples of application of the \texorpdfstring{$d$}{d}-transformation
}

In many problems of applied mathematics,  
the solution is expressed as an infinite series  
\[
\sum_{r=0}^\infty a(r) \varphi_r(x),
\]
where \( \varphi_r(x) \) are orthogonal polynomials,  
elementary functions, special functions,  
or products of such functions.  

Typically, the functions \( \varphi_r \) satisfy  
a linear recurrence relation of finite order,  
and the sequence \( \{\varphi_r(x)\} \)  
belongs to a class \( B^{(m)} \) for some \( m \).  
As follows from Lemma \ref{LemmaB},
if \( a(r) = r^\nu \) or \( a \in {A}^{(\gamma)}\) for some \( \gamma \),  
then, the sequence \( \{ a(r) \varphi_r(x) \} \)  
also belongs to \( {B}^{(m)} \).

We now illustrate the application  
of the \( d \)-transformation to several  
infinite series of the type described  
in the preceding paragraph.  
In particular, we consider series  
whose associated sequences lie in  
\( B^{(2)} \) and \( B^{(4)} \). We note that when the sequence lies in \( B^{(1)} \),  
the \( d \)-transformation reduces to  
Levin's \( u \)-transformation.

\begin{example}\label{Legendre} 
Consider the series  
\[
g(x)=\sum_{n=0}^\infty \frac{P_n(x)}{(1 - 2n)(2n + 3)} = \sqrt{\frac{1 - x}{8}},
\]  
where \( P_n(x) \) denotes the Legendre polynomial of degree \( n \).  
This series converges slowly for \( |x| \leq 1 \), and diverges for \( |x| > 1 \).  

The Legendre polynomials satisfy the three-term recurrence relation 
\begin{equation}\label{Pn}
(n + 1) P_{n+1}(x) = (2n + 1) x P_n(x) - n P_{n-1}(x).
\end{equation}
Using the recurrence relation~\eqref{Pn}, it follows that the sequence \( \{P_n(x)\} \) belongs to \( B^{(2)} \).  
Define
\[
f(n) = \frac{P_n(x)}{(1 - 2n)(2n + 3)}.
\]
The sequence \( 1/(1 - 2n)/(2n + 3)\), being a rational function of \( n \), belongs to \( B^{(1)} \).  
Therefore, by Lemma~\ref{LemmaB}, the sequence \( \{f(n)\} \) 
belongs to the class \( B^{(2)} \). This provides all the 
necessary information to apply the \( d \)-transformation, 
specifically the \( d^{(2)} \)-transformation in this case.

Table~\ref{Table1} shows the results of applying 
the \( d^{(2)} \)-transformation. Recall that 
\( d_r^{(2)} \) requires only \( 2r + 1 \) terms of the series.

The convergence is less impressive at \( x = 0.9 \), 
which is close to the branch point of \( g \). 
In contrast, rapid convergence is observed at \( x = -1.5 \), 
a point where the original series diverges.

\begin{table}[ht]

\centering
\caption{Values of the approximations \( d_r^{(2)}(\{A_N(f)\}) \) to \( g(x)\).}
\medskip
\begin{tabular}{|c|c|c|c|}
\toprule
\( r \) & \( d_r^{(2)},\ \ x=-1.5 \) & \( d_r^{(2)},\ \ 0=0.5 \) & \( d_r^{(2)},\ \ x=0.9 \) \\
\midrule
2   & 0.559015         & 0.2505         & 0.116 \\
4   & 0.559016998      & 0.249998       & 0.1114 \\
6   & 0.5590109944372  & 0.24999989     & 0.11177 \\
8   & 0.55901699437493 & 0.2499999978   & 0.111800 \\
10  & 0.55901099437485 & 0.250000000027 & 0.1118039 \\
\midrule
Exact & 0.559016994374947 & 0.25 & 0.111803398874 \\
\bottomrule
\end{tabular}
\label{Table1}
\end{table}

\end{example}

\begin{example}\label{Ex2}

Another interesting example from \cite{LevinSidi1981} is the following:

Consider the infinite series

\begin{equation}\label{Eq22}
q(\beta,\phi)=\sum_{n=0}^\infty \cos((n+1/2)\beta)P_n(cos\phi)=\begin{cases}
  1/\sqrt{2(cos\beta-cos\phi)} & , \ 0\le\beta<\phi<\pi \\
  \ \ \ \ \ \ \ \ \ \ \ 0  & ,\ 0<\phi<\beta<\pi.
\end{cases}
\end{equation}
\end{example}

Such a series arises in the analysis of solutions to Laplace’s equation on the sphere, under the assumption of spherical symmetry, using separation of variables. The series converges very slowly, and to accelerate its convergence, one may apply the 
$d$-transformation. The first step in this process is to identify the class of functions to which the terms of the series belong. The sequence 
$
\left\{ \cos\left( \left(n + \tfrac{1}{2} \right)\beta \right) \right\}
$
satisfies a linear three-term recurrence relation with constant coefficients, 
and thus belongs to the class \( B^{(2)} \). 
As discussed above, the sequence of Legendre polynomials 
$
\left\{ P_n(\cos\phi) \right\}
$
also belongs to $ B^{(2)}$. 
Therefore, for any fixed \( \beta \) and \( \phi \), the product sequence 
\[
f_n = \cos\left( \left(n + \tfrac{1}{2} \right)\beta \right) P_n(\cos\phi)
\]
belongs to the class \( B^{(4)} \), by Lemma~\ref{LemmaB}. Indeed, as shown in Table~\ref{Table2}, applying the \( d^{(4)} \)-transformation
to the series \eqref{Eq22} yields rapid convergence
to the values of the target function.

The table presents approximations computed 
for various values of \( r \), demonstrating a substantial 
acceleration in convergence and high accuracy 
even for small \( r \). In particular, with \( r = 6 \), 
the value \( d^{(4)}_6 \) utilizes only 25 terms of the series 
and already yields machine-precision approximations 
to \( q(\beta, \phi) \).

\begin{table}[ht]
\centering
\caption{Values of the approximations 
\( d_r^{(4)}\) to \(q(\beta,\phi)\)
for \( \beta = \tfrac{2\pi}{3}, \phi = \tfrac{\pi}{6} \) 
and \( \beta = \tfrac{\pi}{6}, \phi = \tfrac{2\pi}{3} \).}
\medskip
\begin{tabular}{| c|c|c |}
\toprule
$r$ & $d_r^{(4)},\ \beta= 2\pi/3,\ \phi=\pi/6$ & $d_r^{(4)},\ \beta=\pi/6,\ \phi= 2\pi/3$ \\
\midrule
2 & $\ \ \sim 4 \times 10^{-6}$ & 0.604998 \\
3 & $\sim-2\times 10^{-8}$ & 0.60500026 \\
4 & $\sim-2 \times 10^{-10}$ & 0.60500033358 \\
5 & $\sim-9 \times 10^{-13}$ & 0.6050003337080 \\
6 & $\sim-2 \times 10^{-14}$ & 0.605000333706045 \\
\midrule
Exact & $0$ & 0.605000333706055 \\
\bottomrule
\end{tabular}
\label{Table2}
\end{table}

\section{The \texorpdfstring{$D$}{D}-transformation for infinite integrals}

We now present the derivation of the \( D \)-transformation, 
which is designed to accelerate the convergence of infinite integrals 
of the form \( I = \int_{0}^\infty f(t) dt \). 
As in the case of infinite series, the construction of the transformation 
is based on an appropriate asymptotic expansion for the remainder term 
of the integral. Specifically, we consider
\begin{equation}\label{RIN}
R_I(x) = I - \int_x^\infty f(t)dt,
\end{equation}
where \( R_I(x) \) denotes the remainder after truncating the integral at \( x \).

In analogy with the class \( B^{(m)} \) 
defined for infinite series, 
we now introduce a corresponding class 
of functions suitable for applying 
the \( D \)-transformation to infinite integrals.

\begin{definition}{\bf The set $\tilde B^{(m)}$}

Let \( m \in \mathbb{N} \). The set \( \tilde B^{(m)} \) consists of all functions $f$ on $[0,infty)$ which satisfy a linear differential equation of order \( m \), that is,
\begin{equation}\label{ode}
f(t) = \sum_{k=1}^{m} p_k(t)f^{(k)}(t),
\end{equation}
where the coefficient functions \( p_k(n) \) belong to the asymptotic class \( A^{(k)} \), for \( k = 1, 2, \dots, m \). We assume that $m$ is minimal.
\end{definition}

\begin{remark}
We note that most of the special functions used in applied mathematics are in $\tilde B^{(2)}$.
\end{remark}

In analogy to Lemma \ref{LemmaB}, the following lemma, also proved in \cite{LevinSidi1981}, assists in identifying the set $\tilde B^{(m)}$ to which certain functions $f$ belong.

\vspace{1em}

\begin{lemma}\label{LemmatB}{\bf Algebraic rules for the classes \(\{\tilde B^{(m)}\)\} }

If $f\in \tilde B^{(m)}$ and $g\in \tilde B^{(k)}$ then
$f+g\in \tilde B^{(m+k)}$,
$f\cdot g\in \tilde B^{(mk)}$ and $f^2\in \tilde B^{(\frac{m(m+1)}{2})}$.
\end{lemma}

The following theorem, proved in \cite{LevinSidi1981}, is the principal result underlying the definition of the 
$D$-transformation. It establishes an asymptotic expansion for the remainder of an infinite integral whose term sequences belong to the class  $\tilde B^{(m)}$.

\begin{theorem}\label{TRI}
Let \( f \) be integrable on \( [0, \infty) \) and let $f\in \tilde{B}^{(m)}$. Assume that
\[
\lim_{x \to \infty} p_k^{(i)}(x) f^{(k-i)}(x) = 0, \quad \text{for } k = i, i+1, \dots, m, \quad i = 1, 2, \dots, m,
\]
and
\begin{equation}
\sum_{k=1}^m\ell(\ell-1)\cdots(\ell-k+1)p_{k,0}\ne 1,\ \ell\ge -1, \ \ell\ integer,
\end{equation}
where \( p_{k,0} = \lim_{x \to \infty} x^{-i_k} p_k(x) \). Then, as \( x \to \infty \), the remainder $R_I(x)$ satisfies the asymptotic relation
\begin{equation}\label{RemI}
R_I(x)=
\int_x^\infty f(t) \, dt \sim \sum_{k=0}^{m-1}
f^{(k)}(x)
\sum_{j=0}^{\infty} \beta_{k,j} x^{k - j}.
\end{equation}
\end{theorem}

\begin{remark}\label{Remarkp}
The proof of Theorem~\ref{TRI} in~\cite{LevinSidi1981}  
is based on integration by parts of \( R_I(x) \),  
using the differential equation model~\eqref{ode}  
satisfied by \( f \). 
Similar to the proof of Theorem~\ref{RemI},  
the proof of Theorem~\ref{TRS} in~\cite{LevinSidi1981}  
uses ``summation by parts'' applied to the remainder \( R_S(N) \).  
In~\cite{LevinSidi1981}, the discussion begins with infinite integrals,  
as the analysis based on integration by parts is simpler in that case.
\end{remark}

\subsection{Defining the \texorpdfstring{$D$}{D}-transformation for infinite integrals}

We now us the asymptotic expansion~\eqref{RemI} to derive the \( D \)-transformation for infinite integrals $I = \int_{0}^\infty f(t)$dt. We truncate each of the \( m \) infinite sums appearing in~\eqref{RemI} after \( r \) terms. Recalling that \( R_I(x) = S - \int_{0}^{x} f(t)dt \) and denoting the partial integrals of the infinite integral by $A_x(f)=\int_{0}^{x} f(t)dt$, we obtain, up to the truncation error,
\begin{equation} \label{RemTI}
I \approx A_x(f)+ \sum_{k=0}^{m-1}
f^{(k)}(x)
\sum_{j=0}^{r-1} \beta_{k,j} x^{k - j}..
\end{equation}
There are \( mr + 1 \) unknown parameters in~\eqref{RemTI}, namely the integral \( I \) and the coefficients \( \{ \beta_{k,i} \} \).
To determine these parameters, we construct a system of \( mr + 1 \) equations of the form
\begin{equation}
\label{RemTIxj}
\bar I = A_{x_j}(f) + \sum_{k=0}^{m-1}
f^{(k)}(x_j)
\sum_{i=0}^{r-1} \bar\beta_{k,i} x_j^{k - i}
\quad j = 1,2,\ldots, mr+1,
\end{equation}
where \( 0 < x_1 < x_2 < \cdots < x_{mr+1} \). A common choice of evaluation nodes is \( x_j = \ell + jh \),  
for \( j = 1,2,\ldots, mr+1 \). 
With this choice, we solve system~\eqref{RemTNj},
and denote the result for \( \bar{I} \) by \( D^{(m)}_{r,\ell,h} \).

This quantity is known as the \( D \)-\emph{transformation}  
of the sequence of partial integrals of the infinite integral,  
and is denoted by \( D^{(m)}_{r,\ell,h}(\{A_x(f)\}) \).  
The parameters involved are:  
\begin{itemize}
  \item \( m \): the order of the difference operator;  
  \item \( r \): the number of terms used in the asymptotic expansions;  
  \item \( \ell,h \): parameters determining the sequence of evaluation nodes $X=\{x_j\}_{j=1}^{mr+1}$.
  .
\end{itemize}

Note that computing \( D^{(m)}_{r,\ell.h} \)  
uses values of $f$ on the interval $[0,\ell+(mr+1)h]$.  
For the special case \( \ell = 0, h=1 \),  
we write the transformation as  
\( D^{(m)}_{r}(\{A_x(f)\}) \).

For other choices of evaluation nodes $X=\{x_j\}_{j=1}^{mr+1}$, we denote the transformation by $D^{(m)}_{r,X}(\{A_x(f)\})$.

\subsection{Examples of application of the \texorpdfstring{$D$}{D}-transformation}

In many applied mathematics problems, one encounters the need to compute infinite integrals involving special functions.

\begin{itemize}
  \item \textbf{Integral Transforms:}  
  Arise in Fourier and Laplace analysis.  
  Involve exponential, sine, and gamma functions.

  \item \textbf{Green’s Functions:}  
  Used in boundary value problems.  
  Often involve Bessel or Legendre functions.

  \item \textbf{Quantum Mechanics:}  
  Scattering and wavefunction integrals.  
  Include Bessel, Airy, and exponential terms.

  \item \textbf{Statistical Mechanics:}  
  Partition functions and density of states.  
  Often expressed using gamma-type integrals.

  \item \textbf{Optics and Diffraction:}  
  Fresnel and sinc-type integrals.  
  Model light propagation and interference.

  \item \textbf{Fluid Dynamics:}  
  Velocity and pressure fields.  
  Include integrals of Bessel functions.

  \item \textbf{Probability and Statistics:}  
  Moments and distribution tails.  
  Use gamma, beta, and error functions.

  \item \textbf{Asymptotic Analysis:}  
  Study of behavior at infinity.  
  Involves log, sine, and rational kernels.
\end{itemize}

We now illustrate the application  
of the \( D \)-transformation to several  
types of infinite integrals.  
In particular, we consider integrals 
whose associated integrands lie in  
\( \tilde B^{(2)} \) and \( \tilde B^{(3)} \).

\begin{example}\label{Ex3}

An interesting example from \cite{LevinSidi1981} involves the class of the following highly oscillatory infinite integrals:
\begin{equation}
I(a,b)=\int_0^\infty sin(at^2+bt)dt.
\end{equation}
Denoting $f(t)=sin(at^2+bt)$ it turns out that $f$ satisfies the second-order differential equation
$$f(t)=[2a/(2at+b)^3]f'(t)-[1/(2at+b)^2]f''(t).$$
Therefore, $f\in \tilde B^{(2)}$, and we apply the $D^{(2)}$ transformation. In Table \ref{tableab} we present results for the two cases, $a=\pi/2,\ b=0$ and $a=\pi/2,\ b=\pi/2$, in view of the exact values $I(\pi/2,0)=0.5$ and $I(\pi/2,\pi/2)=0.3992050585256$.
\begin{table}[ht]
\centering
\caption{Values of the Approximations \( D_{r,0.2,0.2}^{(2)}\)\ , where \( I(a,b) = \int_0^x \sin(a t^2 + b t) \, dt \).}
\medskip
\begin{tabular}{|c|c|c|}
\hline
\textbf{r} & \( a=\pi/2,\ b=0  \) & \( a=\pi/2,\ b=\pi/2  \) \\
\hline
2  & 0.12         & 0.46            \\
4  & 0.495        & 0.397           \\
6  & 0.4993       & 0.399212        \\
8  & 0.5\,00001   & 0.399205044     \\
10 & 0.49999999989 & 0.399205058518  \\
\hline
\end{tabular}
\label{tableab}
\end{table}
Using values of the integrand in the interval \( [0, 4.4] \),  
the \( D \)-transformation achieves  
an accuracy of 11 correct digits. Due to the non-uniform oscillations of \( f \),  
previous transformations for infinite series,  
such as the \( G \)-transformation in~\cite{Gray},  
fail to be effective in this case.

\end{example}

\begin{example}\label{Ex4}

Another example from~\cite{LevinSidi1981} involves  
an infinite integral of a product of Bessel functions:  
\( f(t) = {J_0(t) J_1(t)}/{t} \).

Since \( J_0 \) and \( J_1 \) belong to \( \tilde B^{(2)} \),  
and \( 1/t \in \tilde B^{(1)} \),  
Lemma~\ref{LemmatB} implies that  
\( f \in \tilde B^{(4)} \).

However, it can be shown that \( f \in \tilde B^{(3)} \),  
and we therefore apply the \( D^{(3)}_r \)-transformation. The integrand decays like \( t^{-1.5} \),  
yet the \( D^{(3)} \)-transformation achieves  
11-digit accuracy using function values 
in the interval \( [0, 32] \).

\begin{table}[ht]
\centering
\caption{The approximations \( D_r^{(3)} \) to \( I = \int_0^\infty {J_0(t)J_1(t)}/{t}dt=2/\pi=0.63661977236758 \).}
\medskip
\begin{tabular}{|c|c|}
\hline
\textbf{r} & \( D_r^{(3)} \) \\
\hline
2  & 0.6341         \\
4  & 0.6366097      \\
6  & 0.63661991     \\
8  & 0.63661977204  \\
10 & 0.636619772340 \\
\hline
\end{tabular}
\end{table}

\end{example}

\begin{example}\label{Ex5}

The third example, taken from~\cite{LevinSidi1981},  
concerns an infinite integral of a slowly decaying monotone function:  
\[
f(t) = \frac{\log(1 + t)}{1 + t^2}.
\]

The function \( f \) is integrable at infinity  
and satisfies the differential equation  
\[
f(t) = -\frac{5t^2 + 4t + 1}{2(2t + 1)} f'(t)  
       - \frac{(t^2 + 1)(t + 1)}{2(2t + 1)} f''(t).
\]

Therefore, \( f \in \tilde{B}^{(2)} \),  
and we apply the \( D^{(2)} \)-transformation.  

We emphasize again that the specific form  
of the differential equation is not used—  
only its order is relevant to the application.

\begin{table}[ht]
\centering
\caption{Values of the approximations \( D^{(2)}_{r,X} \)}
\begin{tabular}{|c|c|}
\hline
\textbf{r} & \( D_{r,X}^{(2)}\) \\
\hline
2  & 1.14       \\
4  & 1.46085        \\
6  & 1.46042     \\
8  & 1.46036208   \\
10 & 1.4603621191   \\
\hline
\end{tabular}
\end{table}

\medskip
 The exact value of the integral is
$\approx 1.460362116753$. The approximations were obtained using 
\(X=\{ x_j = e^{0.2(j-1)}\}_{j=1}^{2r+1} \), with \( r = 2, 4, 6, 8, 10 \).

\end{example}

\section{Conclusion}

The work of Levin and Sidi represents a foundational contribution  
to both the theory and practice of convergence acceleration.  
Their nonlinear transformations significantly broaden the class of problems  
that can be treated efficiently, and remain highly relevant  
in contemporary numerical computation.

In all tested cases, the \( d \)- and \( D \)-transformations  
consistently outperform traditional methods,  
even when the latter are specifically adapted to the problem at hand.  
Moreover, these transformations often succeed in cases  
where other methods—such as the \( G \)-transformation  
or the \( \varepsilon \)-algorithm—fail.

For comprehensive accounts of the 
\( d \)- and \( D \)-transformations, 
including their extension to the 
acceleration of vector sequences, see 
the works of Sidi 
\cite{Sidi1982,Sidi1988,Sidi1990,%
Sidi1997,Sidi2003,Sidi2017}, and for 
extensions to double series and integrals, 
see \cite{GreifLevin}.

\end{document}